\documentclass{statsoc} 

\usepackage{amssymb}
\renewcommand\theenumi{\arabic{enumi}}
\title[Marginalization paradox]{The marginalization paradox does not \\
imply inconsistency for improper priors}

\author[T. C. Wallstrom]{Timothy C. Wallstrom} 

\coaddress{Timothy C.\ Wallstrom, 
Theoretical Division, MS-B213, Los Alamos National
Laboratory, Los Alamos, New Mexico, 87545, USA} 

\email{tcw@lanl.gov}

\address{Los Alamos National Laboratory, Los Alamos, New Mexico, USA}

\begin{document}

\begin{abstract}  
  
  The marginalization paradox involves a disagreement between two
  Bayesians who use two different procedures for calculating a
  posterior in the presence of an improper prior. We show that the
  argument used to justify the procedure of one of the Bayesians is
  inapplicable. There is therefore no reason to expect agreement, no
  paradox, and no evidence that improper priors are inherently
  inconsistent. We show further that the procedure in question can be
  interpreted as the cancellation of infinities in the formal
  posterior.  We suggest that the implicit use of this formal
  procedure is the source of the observed disagreement.
  
\end{abstract}

\section{Introduction.}
\label{sec:intro}

An important question in statistics is whether Bayesian inference can
be extended to the setting of improper priors in a consistent and
intuitively viable manner. The use of improper priors was common
throughout much of the twentieth century, and appears to be a useful
idealization for many applications. In the 1970s, however, two
influential arguments appeared against the use of improper priors: the
``marginalization paradox,'' and ``strong inconsistency.''  These
arguments appear to have convinced most statisticians that improper
priors must be abandoned.

In this paper we discuss the marginalization paradox, due
to~\citet*{DSZ73}~(DSZ73).  Let $p(x|\theta)$ be a normalized sampling
distribution with parameter $\theta=(\eta, \zeta)$ and data $x=(y,z)$,
and let $p(\theta)$ be a prior, which may be improper, \textit{i.e.},
of infinite total probability. The marginalization paradox concerns
the problem of calculating $p(\zeta|z)$, under a certain set of
assumptions. A first Bayesian, $B_1$, eliminates $\eta$ and then $y$;
a second Bayesian, $B_2$, eliminates $y$ and then $\eta$. The details
of the procedures are given in DSZ73. It is claimed that these
procedures rely only on principles that would have to hold in any
intuitively viable theory of inference.  If $p(\theta)$ is improper,
however, $B_1$ and $B_2$ generally get incompatible answers. It has
been widely inferred that any extension of Bayesian inference to the
context of improper priors will be inconsistent.

The purpose of this paper is to show that the marginalization paradox
does not imply that the use of improper priors will lead to
inconsistency. First, we show that the argument used to justify
$B_1$'s elimination of $y$ is invalid, because it is based on the
application of probabilistic intuitions to a formal quantity whose
probabilistic meaning has not been justified. The ``paradox'' is
thereby resolved, since we now have no reason to believe that $B_1$'s
answer is correct, and no reason to insist that the answers of $B_1$
and $B_2$ be compatible.

Second, we analyze $B_1$'s procedure on its own terms, to get a better
sense for what is being assumed. The posterior $p(\zeta|z)$ is defined
as a ratio, which is only formal when the prior is improper because
there are infinities in the numerator and denominator. $B_1$'s
procedure is equivalent to the assumption that these infinities will
cancel. What DSZ73 have shown, therefore, is that there is no
consistent extension of Bayesian inference in which the cancellation
law, assumed implicitly by $B_1$, holds when the prior is improper.
But this is only to be expected: it is analogous to the well-known
fact that there is no consistent extension of arithmetic to the
extended real numbers in which the cancellation law holds for
infinity. The proposal that we abandon improper priors because of the
marginalization paradox is analogous to the proposal that we abandon
the use of infinity because it does not obey the laws of arithmetic.

In brief, the inconsistency of the marginalization paradox is based on
an assumption that has not been justified intuitively and that is
unreasonable mathematically. There is nothing in the marginalization
paradox to preclude the existence of a formalism that justifies the
careful use of improper priors.

\section{The intuitive argument.}
\label{sec:ia}

In this section we show that the validity of $B_1$'s argument has not
been established, because it is based on an intuitive probabilistic
argument, and the distribution to which it is applied has not been
shown to have a probabilistic meaning.  In other words, we show that
DSZ73 have not made their case, because their argument contains a gap.

In addition to the assumptions described in Section~\ref{sec:intro},
we assume the following:
\begin{enumerate}
\item The formal posterior, defined as
\begin{equation}
    \label{eq:post}
    p(\zeta|y,z) = {\int p(y,z|\eta,\zeta) \,p(\eta,\zeta)
                 \,d\eta\over \int p(y,z|\eta,\zeta) \,p(\eta,\zeta)
                 \,d\eta\,d\zeta},
\end{equation}
is independent of $y$. We denote the common value by $p_1(\zeta|z)$.
Note that the value of $p(\zeta|y,z)$ and the validity of the
assumption itself depend on the prior. \bigskip

\item  The marginalized sampling distribution,
   \begin{displaymath}
    p(z|\eta,\zeta) = \int p(y,z|\eta,\zeta)\,dy
  \end{displaymath}
  is independent of $\eta$. We denote the common value by
  $p_2(z|\zeta)$. \bigskip
  
\item For each value of $\zeta$, the prior is improper in $\eta$: $\int
  p(\eta,\zeta) \, d\eta = \infty$.

\end{enumerate}
Assumptions 1 and 2 enable $B_1$ and $B_2$, respectively, to invoke
intuitive arguments to determine $p(\zeta|z)$, even though the formal
calculations would lead to infinities. Assumption~3 is
satisfied by all of the examples in DSZ73, and reflects the fact
that we are really interested in impropriety in $\eta$.

We focus on only one aspect of the analysis in DSZ73, because we
believe that aspect to be the source of all of the difficulties. The
aspect in question is $B_1$'s elimination of $y$, which occurs after
he has already marginalized over $\eta$.  $B_1$ assumes that since
$p(\zeta|y,z)$ is independent of $y$, then $p(\zeta|z)$ must be equal
to the $y$-independent value of this function. 
  
The justification that DSZ73 give for this assumption is intuitive,
and has been formalized as the ``reduction principle,'' which is
stated as follows in~\cite*{DSZ96}: ``Suppose that a general method of
inference, applied to data $(y,z)$, leads to an answer that in fact
depends on $z$ alone.  Then we should obtain the same answer if we
apply the method to $z$ alone.'' The principle enables one to
determine the answer to the problem with data $z$ from the answer to
the problem with data $(y,z)$, provided that the latter answer depends
only on $z$.  We have no objection to this principle as stated.  We
wish to emphasize, however, that in order to apply the principle (or
invoke the intuition behind the principle), we must first have the
``answer'' to a problem of inference, given data $(y,z)$.

The problem with $B_1$'s argument is that $p(\zeta|y,z)$ has not been
shown to be the ``answer'' to a problem of inference, so the reduction
principle is inapplicable. We show below that in the context of the
marginalization paradox, any sampling distribution $p(y,z|\zeta)$
associated with $p(\zeta|y,z)$ is necessarily \textit{improper}, so
that it has no inherent probabilistic meaning.  There is no reason to
assume that the associated formal posterior will have any
probabilistic meaning, even if that posterior is proper.  In the
absence of such a meaning, $p(\zeta|y,z)$ is not the answer to a
problem of inference, $B_1$ is unable to use the reduction principle
to complete his argument, and the inconsistency vanishes.

We are not claiming that it is impossible to provide a meaning for an
improper distribution. Indeed, such an assumption would preclude the
use of improper priors and prejudge the whole issue. We are merely
observing that in order to use the reduction principle, a
probabilistic meaning must be provided for $p(y,z|\zeta)$, and this
has not been done.  Even if a meaning is provided, any manipulations
of the distribution must be justified in terms of that meaning, and
there is no guarantee that the resulting procedures will be the formal
analogs of valid procedures for proper distributions.

We now establish the impropriety of the sampling distribution.  

\noindent\textbf{Proposition:\ \ } 
Let $p(\eta,\zeta)$ be given, and let $p(\eta,\zeta)= p(\eta|\zeta)\,
p(\zeta) $ be any factorization of $p(\eta,\zeta)$ such that
$0<p(\zeta)<\infty$. Under the above assumptions we have, for
each $\zeta$,
\begin{equation}
  \int p(y,z|\zeta)\, dy = \infty.
    \label{eq:result}
\end{equation}

\noindent\textit{Proof:\ \ } 
\begin{eqnarray*}
  \int p(y,z|\zeta) \,dy = {1\over p(\zeta)} \int p(y,z|\eta,\zeta) 
    p(\eta,\zeta)\,d\eta
   \,dy = {p_2(z|\zeta)\over p(\zeta)} 
   \int p(\eta,\zeta) \,d\eta =    \infty. \qquad
\end{eqnarray*}
The interchange in the order of integration is justified by
Tonelli's theorem. $\square$

An immediate corollary is that $\int p(y,z|\zeta)\, dy\,dz = \infty$.
The factorization of $p(\eta,\zeta)$ is nonunique, and this implies a
nonuniqueness in the definition of $p(y,z|\zeta)$. The proposition
shows, however, that impropriety of the conditional distribution is
independent of the choice of factorization. Note also that although we
are evaluating $B_1$'s argument, the proof depends on assumption~(2),
which was made for $B_2$'s benefit.

\section{The formal argument.}
\label{sec:fa}

We now consider $B_1$'s procedure on its own terms, as a formal
procedure. We find that in the case of a proper prior, $B_1$'s use of
the reduction principle is equivalent to the cancellation of a finite
factor in a ratio defining $p(\zeta|z)$, and in the case of an
improper prior, to the cancellation of an infinite factor.  It is
well-known that the formal cancellation of infinities will generally
lead to inconsistencies. We conclude that when viewed formally,
$B_1$'s procedure is highly suspect.

In general, the posteriors of $\zeta$ given $(y,z)$ and given $z$ are
given formally by the following expressions:
\begin{eqnarray}
  \label{eq:postx}
  p(\zeta|y,z) &=& {p(y,z,\zeta)\over \int p(y,z,\zeta) \,d\zeta},
    \qquad\qquad \mathrm{and} \\
  \label{eq:postz}
  p(\zeta|z) &=& {\int p(y,z,\zeta)\,dy\over \int p(y,z,\zeta)
  \,dy\,d\zeta}.
\end{eqnarray}
Under Assumption~1, $p(\zeta|y,z)$ is independent of $y$. Then
\begin{equation}
  \label{eq:fact}
  p(y,z,\zeta) = p(y,z)\, p_1(\zeta|z),
\end{equation}
where $p(y,z) =\int p(y,z,\zeta) d\zeta$.  
Substituting Eq.~(\ref{eq:fact}) into Eq.~(\ref{eq:postz}), we obtain
\begin{equation}
  p(\zeta|z) = {\int p(y,z) \,p_1(\zeta|z)\, dy \over \int 
         p(y,z) \,p_1(\zeta|z) \,dy\,d\zeta}.
\end{equation}
When $\int p(y,z) dy$ is finite, then $p(\zeta|z) = p_1(\zeta|z)$.  

If we also make Assumptions~2 and~3, the proposition implies that
$\int p(y,z)\,dy = \infty$. The assumption that
$p(\zeta|z)=p_1(\zeta|z)$ is now equivalent, as claimed, to the
assumption that it is permissible to cancel infinite factors of
$\int p(y,z)\, dy$ from the ratio defining $p(\zeta|z)$.

\section{Discussion.}

We have observed that the inconsistencies uncovered in DSZ73 depend on
formal manipulation on the part of $B_1$. We have shown, in
Sections~\ref{sec:ia} and~\ref{sec:fa}, respectively, that $B_1$'s
procedure has not been justified intuitively, and is suspect
mathematically. We therefore see no reason to accept $B_1$'s
reasoning, or to regard the validity of this reasoning as necessary or
desirable in any extension of Bayesian inference to improper priors.
Once $B_1$'s reasoning is rejected, the marginalization paradox
disappears.

The core of our argument is the observation that $B_1$'s argument is
formal because the sampling distribution $p(y,z|\zeta)$ is improper.
To the best of our knowledge, this observation has not been made
previously. The impropriety of the sampling distribution has perhaps
been obscured by its nonuniqueness and by the fact that the formal
posterior can be calculated from Eq.~(\ref{eq:post}) without ever
computing the sampling distribution explicitly.

Previous analyses of the marginalization paradox generally accepted
the validity of both Bayesians' arguments. The problem then becomes
one of understanding when and why the two Bayesians will agree. This
analysis was initiated in DSZ73, which is mostly dedicated to this
question.  It turns out that for problems amenable to group analysis,
consistency may be achieved by a uniquely determined prior.  The
priors determined by this constraint, however, are unsatisfactory for
a variety of reasons, which DSZ73 explore in detail.  They conclude
that an acceptable theory is elusive or unachievable.

The most persistent and insightful critic of the marginalization
paradox has been the late E.~T.~Jaynes.
Cf.~\cite{Jaynes80a,DSZ80,Jaynes80b,DSZ96,Jaynes03}, for his extended
debate with the authors of DSZ73.  We believe that at the conceptual
level, Jaynes' critique was fundamentally correct, in that he
identified the source of the inconsistencies as the formal
manipulation of completed infinities.  A particularly elegant
statement of this view can be found in~\cite{Jaynes03}.  At the
technical level, Jaynes did not recognize that $B_1$'s argument was
invalid, so he was forced to try to determine how the two Bayesians
could be reconciled. His thesis was that the disagreement between the
Bayesians reflected differences in their prior information. In our
opinion, this analysis was not entirely successful, and the correct
approach is to reject $B_1$'s reasoning.

For general background on the marginalization paradox and related
issues, we refer the reader to the excellent review article
of~\cite{KW96}.

\section*{Acknowledgments}  I thank Harry Martz, Brad Plohr, and
Arnold Zellner for helpful comments. I also acknowledge support from
the Department of Energy under contract W-7405-ENG-36.

\end{document}